\documentclass[11pt, twoside]{article}
 \usepackage{amsfonts, amssymb, amsmath, amsthm, latexsym, amscd}
\newcommand{\braket}[2]{\langle #1,#2 \rangle}

\newcommand{\dDel}{\bigtriangledown}

\newcommand{\Ga}{\Gamma}
\newcommand{\ga}{\gamma}
\newcommand{\Om}{\Omega}
\newcommand{\del}{\delta}
\newcommand{\Del}{\Delta}
\newcommand{\La}{\Lambda}

\newcommand{\tht}{\theta}
\newcommand{\eps}{\epsilon}

\newcommand{\f}{\frac}

\newcommand{\lo}{\longrightarrow}

\newcommand{\BR}{\mathbb{R}}
\newcommand{\al}{\alpha}

\newcommand{\pa}{\partial}

\theoremstyle{plain}
\newtheorem{Th}{Theorem}[section]
\newtheorem{Lem}[Th]{Lemma}
\newtheorem{Prop}[Th]{Proposition}

\theoremstyle{definition}

\begin{document}
\addtolength{\textheight}{0 cm} \addtolength{\hoffset}{0 cm}
\addtolength{\textwidth}{0 cm} \addtolength{\voffset}{0 cm}

\title{Soliton solutions for  quasilinear Schr\"{o}dinger
equations involving  supercritical   exponent in ${\BR}^N$}

\author{Abbas Moameni  \textsc \bf \footnote{Research is supported by a
Postdoctoral Fellowship at the University of British Columbia.}}
\maketitle
\begin{center}
{\small Department of Mathematics \\
\small  University of British Columbia\\
\small Vancouver, B.C., Canada \\
\small {\tt moameni@math.ubc.ca }}
\end{center}

\begin{abstract}
We study the existence of positive solutions to the quasilinear
elliptic problem
\begin{eqnarray*}
-\epsilon \Delta u+V(x)u-\epsilon k(\Del(|u|^{2}))u=g(u), \quad u>0,
x \in
  {\BR}^N,
\end{eqnarray*}
where g has  superlinear growth at infinity without any restriction
from above on its growth. Mountain pass in a suitable Orlicz space
is employed to establish this result.  These equations contain
strongly singular nonlinearities which include derivatives of the
second order which make the situation more complicated. Such
equations arise when one seeks for standing wave solutions for the
corresponding quasilinear Schr\"{o}dinger equations. Schr\"{o}dinger
equations of this type have been studied as models of several
physical phenomena. The nonlinearity here corresponds to the
superfluid film equation in plasma physics.
\end{abstract}
\section{Introduction}

 We are concerned with the existence of positive
solutions for quasilinear elliptic equations in the entire space,
\begin{eqnarray*}
-\epsilon \Delta u+V(x)u-\epsilon k(\Del(|u|^{2}))u=g(u), \quad u>0,
x \in
  {\BR}^N,
\end{eqnarray*}
where $\epsilon$ is a positive parameter, $V:{\BR}^N\rightarrow
[0,\infty)$ and $g:[0,\infty)\rightarrow [0,\infty)$ are nonnegative
continuous functions. Solutions of this equation are related to the
existence of standing wave solutions for quasilinear Schr\"{o}dinger
equations of the form
\begin{equation}
i\pa_t z=-\epsilon \Del z+W(x)z-l(|z|^2)z-k\epsilon\Del
h(|z|^2)h'(|z|^2)z, \quad x\in {\BR}^N,  N >2,
\end{equation}
where $W(x)$  is a given potential, $k$ is a real constant and $l$
and $h$ are real functions. Quasilinear equations of the form (1)
have been established  in several areas of physics corresponding to
various types of $h$.
 The  superfluid film equation in
plasma physics has this structure   for $h(s)=s$, ( Kurihura in
[8]). In the case $h(s)=(1+s)^{1/2}$, equation (1) models the
self-channeling of a high-power ultra short laser in matter, see
 [21]. Equation (1) also
appears in  fluid mechanics [8,9], in the theory of Heidelberg
ferromagnetism and magnus [10], in dissipative quantum mechanics [7]
and in condensed matter theory [14]. We consider the case  $h(s)=s$
and  $k>0$. Setting $z(t,x)=\exp(-iFt)u(x)$ one obtains a
corresponding equation of elliptic type which has the formal
variational structure:
\begin{eqnarray}
-\epsilon \Delta u+V(x)u-\epsilon k(\Del(|u|^{2}))u=g(u), \quad u>0,
x \in
  {\BR}^N,
\end{eqnarray}
where $V(x)=W(x)-F$ is the new potential function and $g$ is the new nonlinearity.\\

 Note that, for the case $g(u)=|u|^{p-1}u$ with $N\geq 3,$
$p+1=22^*= \frac{4N}{N-2}$
 behaves like a critical exponent for the above equation [13, Remark
 3.13]. For the
 subcritical case $p+1 <22^*$ the existence of solutions for problem (2) was studied in [10, 11, 12, 14, 15,
 16] and it was left open for the critical exponent case $p+1=22^*$ [13; Remark
 3.13].
   The author in [16], proved the
existence of solutions   for $p+1=22^*$ whenever the potential
function $V(x)$ satisfies some geometry conditions. It the present
paper, we will extend this result for the supercritical case. It is
well-known that for the semilinear case $(k=0)$,
\begin{eqnarray*}
-\epsilon \Delta u+V(x)u=g(u), \quad u>0, x \in
  {\BR}^N, \quad \quad \quad
\quad \quad \quad \quad \quad \quad \quad (P)
\end{eqnarray*}

$p+1=2^*$ is the critical exponent when $N \geq 3$. In terms of
critical growth, the case $N=2$ is quite different than $N \geq 3.$
We divide, these studies in three cases for problem ($P$):

\begin{itemize}
 \item {\it  Subcritical growth: } $\lim_{t\rightarrow +\infty} \frac {|g(t)|}{|t|^{2^*}}=0,$ if $N\geq 3$; and
 $\lim_{t\rightarrow +\infty} \frac {|g(t)|}{\exp(\beta t^2)}=0$
 for all $\beta >0,$ if $N=2.$
 \item {\it   Critical growth:  }$\lim_{t\rightarrow +\infty} \frac {|g(t)|}{|t|^{2^*}}=L$ with $L>0$ if $N\geq 3$;
 and for $N=2,$ there exists $\beta_0>0$ such that
 $$\lim_{t\rightarrow +\infty} \frac {|g(t)|}{\exp(\beta t^2)}=0 \quad \quad \forall \beta >\beta_0,
 \quad  \lim_{t\rightarrow +\infty} \frac {|g(t)|}{\exp(\beta t^2)}=+\infty \quad \quad \forall \beta <\beta_0. $$
 \item {\it  Supercritical growth:   }$\lim_{t\rightarrow +\infty} \frac {|g(t)|}{|t|^{2^*}}=+\infty,$ if $N\geq 3$; and
 $\lim_{t\rightarrow +\infty} \frac {|g(t)|}{\exp(\beta t^2)}=+\infty$
 for all $\beta >0,$ if $N=2.$
\end{itemize}
Note that the corresponding critical growth for $N=2$ comes from a
version of Moser-Trudinger inequality in whole space ${\BR}^2$ (see
[6]) as follows,
\begin{eqnarray*}
\int_{{\BR}^2}\big( \exp (\beta |u|^2)-1 \big )\, dx <+\infty, \quad
\quad \quad \forall u \in H^1({\BR}^2), \beta >0.
\end{eqnarray*}
Also, if $\beta< 4\pi$ and $|u|_{L^2({\BR}^2)}\leq C,$ there exists
a constant $C_2=C_2(C,\beta)$ such that
\begin{eqnarray*}
\sup_{|\nabla u|_{L^2({\BR}^2)}\leq 1}\int_{{\BR}^2}\big( \exp
(\beta |u|^2)-1 \big )\, dx <C_2.
\end{eqnarray*}

There are many results about the existence of solutions for the
subcritical, critical  and the supercritical
exponent case for problem ($P$)(e.g. [1, 4, 5, 19, 22]).\\

 In the case $k >0$, for the subcritical case,  the
existence of a nonnegative solution was proved  for $N=1$ by
Poppenberg, Schmitt and Wang in [18] and for $N\geq 2$ by  Liu and
Wang  in [12].
   In [13] Liu and Wang improved
these results by using a change of variables and treating the new
problem in an Orlicz space. The author in [15], using the idea of
the fibrering method, studied this problem in connection with the
corresponding eigenvalue problem for the laplacian $ -\Del u=V(x) u$
and  proved the existence of multiple solutions for  problem (2).
 It is  established in  [11],
the existence of both one-sign and nodal ground states of soliton
type solutions by the Nehari method. They also established some
regularity of the positive solutions.

As it was mentioned above, for the case $k>0$ with $g(u)=|u|^{p-1}u$
and $N \geq 3$, $p+1=22^*= \frac{4N}{N-2}$
 behaves like a critical exponent for problem (2). This is  because of the nonlinearity term  $-\epsilon k(\Del(|u|^{2}))u$.
 Therefore for problem (2), because of the presence of this nonlinearity term,
 the above definition of  Subcritical , Critical  and Supercritical
 growth changes as
   follows:

\begin{itemize}
 \item {\it  \textbf{Subcritical growth: } $\lim_{t\rightarrow +\infty} \frac {|g(t)|}{|t|^{22^*}}=0,$ if $N\geq 3$; and
 $\lim_{t\rightarrow +\infty} \frac {|g(t)|}{\exp(\beta t^4)}=0$
 for all $\beta >0,$ if $N=2.$}
 \item {\it   \textbf{Critical growth: } }$\lim_{t\rightarrow +\infty} \frac {|g(t)|}{|t|^{22^*}}=L$ with $L>0$ if $N\geq 3$;
 and for $N=2,$ there exists $\beta_0>0$ such that
 $$\lim_{t\rightarrow +\infty} \frac {|g(t)|}{\exp(\beta t^4)}=0 \quad \quad \forall \beta >\beta_0,
 \quad  \lim_{t\rightarrow +\infty} \frac {|g(t)|}{\exp(\beta t^4)}=+\infty \quad \quad \forall \beta <\beta_0. $$
 \item {\it \textbf{ Supercritical growth: }  }$\lim_{t\rightarrow +\infty} \frac {|g(t)|}{|t|^{22^*}}=+\infty,$ if $N\geq 3$; and
 $\lim_{t\rightarrow +\infty} \frac {|g(t)|}{\exp(\beta t^4)}=+\infty$
 for all $\beta >0,$ if $N=2.$
\end{itemize}

Here, we shall study problem (2) with $N\geq 2$ and show the
existence of positive solutions when the function $g$ has the
supercritical growth.  Before to state the main result, we fix the
hypotheses on the potential function  $V$ and the function $g.$
Indeed, we assume that the potential function $V$ is radial, that is
$V(x)=V(|x|),$ and satisfies the following conditions:

There exist $0<R_1<r_1<r_2<R_2$ and $\alpha>0$ such that
\begin{enumerate}
\item[{\bf A1:}] $V(x) =0,  \quad \quad \quad \forall x\in
\Omega : =\left\{ x\in {\mathbb{R}}^N : r_1 < |x| < r_2 \right\},$
\item[{\bf A2:}] $V(x) \geq \alpha, \quad \quad \quad   \forall x\in\Lambda^c,$
\end{enumerate}
where $\Lambda=\left\{ x \in {\BR}^N : R_1< |x| < R_2 \right\}.$
Also, we assume  $g$ is continuous and verifies the following
conditions,
\begin{enumerate}
\item[{\bf H1:}] $\lim_{t\rightarrow +\infty} \frac {g(t)}{t}=+\infty.$
\item[{\bf H2:}] The Ambrosetti-Rabinowitz growth condition: There
exists $\theta >2$ such that
\begin{eqnarray*}
0\leq \theta G (t)= \theta \int _0^t g(s) \, ds \leq t g(t), \quad
\quad t \in \mathbb{R}.
\end{eqnarray*}
\item[{\bf H3:}] $\frac {g(t)}{t}$ is non-decreasing with respect to
$t$, for $t>0.$
\item[{\bf H4:}] $\lim_{t\rightarrow 0} \frac {g(t)}{t}=0.$
 \end{enumerate}

Here is our main Theorem.
\begin{Th} Assume Conditions $ H1- H4,  A1$ and $ A2.$ Then,
 there exists $\eps_0>0$, such that for all $\eps\in
(0,\eps_0)$ problem $(2)$ has a nonnegative solution $u_\eps\in
H^1_r({\BR}^N)$ with $u^2_\eps\in H^1_r({\BR}^N)$ and
$$u_\eps(x)\lo 0 \quad as\quad |x|\lo+\infty.$$
\end{Th}
This paper is organized as follows. In Section 2, we reformulate
this problem in an appropriate Orlicz space. In Section 3, we prove
the existence of a solution for   a special deformation of  problem
(2). Theorem 1.1 is proved in Section 4.
\section{Reformulation of the problem and preliminaries}
Denote by $H_{r}^1 ({\BR}^{N})$ the space of radially symmetric
functions in
\begin{eqnarray*}
H^{1,2} ({\BR}^N)= \left \{ u \in L^{2}({\BR}^N) : \dDel u \in
L^{2}({\BR}^N) \right \},
\end{eqnarray*}
and by $C_{0,r}^{\infty} ({\BR}^N)$  the space of radially symmetric
functions in $C_{0}^{\infty} ({\BR}^N).$ Without loss   of
generality, one can assume $k=1$ in problem (2).
 We formally
formulate problem (2) in a variational structure as follows
$$J_\eps(u)=\f{\eps}{2} \int_{{\BR}^N}(1+u^2)| \dDel u|^2dx+\f{1}{2}
\int_{{\BR}^N}V(x)u^2dx- \int_{{\BR}^N} G(u)dx.$$ on the space
$$X=\{u\in H_r^{1,2}({{\BR}}^N): \int_{{{\BR}}^N} V(x) u^2dx<\infty\},$$
which is equipped with the following norm,
$$\|u\|_X=\left \{ \int_{{\BR}^N}| \dDel u|^2dx+ \int_{{{\BR}}^N} V(x) u^2dx\right \}^{\frac{1}{2}}.$$
 Liu and Wang in [13] for the subcritical case,  by making a change of variables treated
this problem in an  Orlicz space. Following their work, we consider
this problem for the supercritical exponent case  in the same Orlicz
space. To  convince the reader we briefly recall some of their
notations and results that are useful in the sequel.

First, we make a change of variables as follows,
$$dv=\sqrt{1+u^2}du, \quad v=h(u)=\f{1}{2}u\sqrt{1+u^2}+\f{1}{2} \ln
(u+\sqrt{1+u^2}).$$ Since $h$ is strictly monotone it has a
well-defined inverse function: $u=f(v)$. Note that
$$h(u)\sim \begin{cases}
u, & |u|\ll 1\\
\f{1}{2}u|u|, & |u|\gg 1, \quad h'(u)=\sqrt{1+u^2},
\end{cases}$$
and
$$f(v)\sim \begin{cases}
v & |v|\ll 1 \\
\sqrt{\f{2}{|v|}}v, & |v|\gg 1, \quad
f'(v)=\f{1}{h'(u)}=\f{1}{\sqrt{1+u^2}}= \f{1}{\sqrt{1+f^2(v)}}.
\end{cases}$$
Also, for some $C_0>0$ it holds
$$L(v):=f(v)^2\sim \begin{cases} v^2& |v|\ll 1,\\
2|v|& |v|\gg 1, \quad L(2v)\leq C_0 L(v),
\end{cases}$$
$L(v)$ is convex, $L'(v)=\f{2f(v)}{\sqrt{1+f(v)^2}}$,
$L''(v)=\f{2}{(1+f(v)^2)^2}>0$.

Using this change of variable, we can rewrite the functional
$J_\eps(u)$ as
$$\bar {J}_\eps(v)=\f{\eps}{2} \int_{{\BR}^N} |\dDel v|^2dx+ \f{1}{2} \int_{{\BR}^N}V(x)
f(v)^2 dx-\int_{{\BR}^N} G(f(v))dx.$$ $\bar {J}_\eps$ is defined on
the space
$$H^1_L ({\BR}^N)=\{v|   v(x)=v(|x|),   \int_{{\BR}^N}|\dDel v|^2dx <\infty, \int_{{\BR}^N}V(x)
L(v)dx<\infty\}.$$ We introduced the Orlicz space (e.g.[20])
$$E_L ({\BR}^N)=\{v| \int_{{\BR}^N}V(x) L(v)dx<\infty\},$$
equipped with the  norm
$$|v|_{E_L ({\BR}^N)}=\inf_{\zeta>0} \zeta(1+\int_{{\BR}^N}(V(x) L(\zeta^{-1}v(x))dx),$$
and define the norm of $H^1_L ({\BR}^N)$ by
$$\|v\|_{H^1_L ({\BR}^N)}= |\dDel v|_{L^2({\BR}^N)}+|v|_{E_L ({\BR}^N)}.$$
Here are  some related facts. See Propositions~(2.1) and (2.2) in
[13] for the proof.
\begin{Prop}
 \begin{enumerate}
 \item[{ (i)}]
 $E_L({\BR}^N)$ is a Banach space.
 \item[{(ii)}]  If $v_n\lo v$ in $E_L({\BR}^N)$, then $\int_{{\BR}^N}V(x)|
L(v_n)-L(v)| dx\lo 0$ and  $\int_{{\BR}^N}V(x)| f(v_n)-f(v)|^2dx\lo
0$.
 \item[{(iii)}]  If $v_n\lo v$ a.e. and $\int_{{\BR}^N}V(x) L(v_n)dx\lo
\int_{{\BR}^N}V(x)L(v)dx$, then $v_n\lo v$ in $E_L({\BR}^N)$.
 \item[{(iv)}]  The dual space $E^*_L({\BR}^N)=L^\infty\cap L_V^2=\{w| w\in L^\infty,
\int_{{\BR}^N}V(x)w^2dx<\infty\}$.
 \item[{ (v)}]  If $v\in E_L({\BR}^N)$, then $w=L'(v)=2f(v)f'(v)\in E^*_L({\BR}^N)$, and
$|w|_{E^*_L}=\sup_{|\phi|_{E_L ({\BR}^N)}\leq 1}(w,\phi)\leq
C_1(1+\int_{{\BR}^N} V(x)L(v)dx)$, where $C_1$ is a constant
independent of $v$.
 \item[{(vi)}] For $N>2$ the map:$v\lo f(v)$ from $H^1_L({\BR}^N)$ into
$L^q({\BR}^N)$ is continuous for $2\leq q\leq 22^*$ and is compact
for $2< q <22^*.$ Also, for $N=2$, this map is compact for $q>2.$
 \end{enumerate}
\end{Prop}

 Hence forth, $\int, H^1, H^1_r,  H^1_L, E_L, L^t, |\cdot|_L$  and $\|\cdot\|$
stand for $\int_{{\BR}^N}$, $H^{1,2}({\BR}^N)$,  $ H^1_r({\BR}^N)$,
$ H^1_L({\BR}^N)$, $ E_L({\BR}^N)$, $L^t({\BR}^N),$ $|\cdot|_{E_L
({\BR}^N)}$ and $\|\cdot\|_{H^1_L({\BR}^N)}$ respectively. In the
following we use $C$ to denote any constant that is independent of
the sequences considered.

\section{Auxiliary Problem}
In this section, we shall show some results needed  to prove Theorem
1.1. Indeed,  we first consider a special deformation $\bar{H}_\eps$
(see (3) in the following) of $\bar{J}_\eps.$ Then, We show that the
functional $\bar{H}_\eps$ satisfies all the properties of the
Mountain Pass Theorem. Consequently,  $\bar{H}_\eps$ has a critical
point for each $\eps>0.$  We shall use this to prove Theorem 1.1 in
the next section. In fact, we will see that  the functionals
$\bar{J}_\eps$ and $\bar{H}_\eps$  will coincide for the small
values of $\eps$. This idea was explored by
Del Pino and Felmer [5].\\

To do this, we shall consider constants $$k>\max \big
\{\f{\tht}{\tht-2},2\big \},  \quad \quad \quad \quad \quad (\tht
\text { is introduced in } H2),$$ and $a$ with $$\f {g(a)}{a}=
\f{\al}{k}, \quad \quad \quad \quad  \quad \quad \quad \quad \quad(
\al \text { is introduced in } A2),$$ and functions
\begin{align*}
\bar{g}(s)&=\begin{cases} g(s), & s\leq a, \\ (\f{\al}{k})s, & s>a,
\end{cases}\\
w(x,s)&=\chi_\La(x)g(s)+(1-\chi_\La(x))\bar{g}(s),
\end{align*}
where $\chi_\La$ denotes the characteristic function of the set
$\La$. Set $W(x,t)=\int_0^t w(x,\zeta)d\zeta$.
 It is easily seen that the function $w$
satisfies the following conditions,
\begin{enumerate}
\item[{\bf G1:}]
$ 0\leq \tht W(x,t) \leq  w(x,t)t, \quad \forall x\in\La, t\geq 0.$
\item[{\bf G2:}] $0\leq 2W(x,t)\leq w(x,t)t\leq \f{1}{k} V(x)t^2,\quad
 \forall x\in {\La}^c, t \in {\BR}.$
\end{enumerate}
Also, it is easy to check that  $w$ satisfies the condition $H2.$ In
the sequel, we denote by $G3,$ the condition $H2$ with $g$ replaced
by $w.$

Now, we study the existence of solutions for the deformed equation,
i.e.
$$-\eps\Del u+V(x)u-\eps(\Del(|u|^2))u = w(x,u), \quad x\in {\BR}^N.$$
which correspond to the critical points of
$$H_{\eps}(u)=\f{\eps}{2}\int (1+u^2)|\dDel u|^2+\f{1}{2} \int V(x)u^2-\int W(x,u)dx.$$
As in Section (2), we can rewrite the functional  $H_{\eps} (u)$ as
a new functional  $\bar H_{\eps} (v)$ with $u=f(v)$ as follows,
\begin{equation}
\bar{H}_{\eps}(v)=\f{\eps}{2} \int |\dDel v|^2 dx+\f{1}{2} \int
V(x)f(v)^2dx-\int W(x,f(v))dx.
\end{equation}
 $\bar {H}_{\eps} (v)$
is defined on the Orlicz space $H^1_L.$ To simplify the writing in
this section, we shall assume $\eps=1, H_1=H$ and $\bar {H}_1=\bar
H.$

The following Proposition states some properties of the functional
$\bar {H}.$
 \begin{Prop}
 \begin{enumerate}
 \item[{(i)}]  $\bar{H}$ is well-defined on $H^1_L$.
 \item[{(ii)}]  $\bar{H}$ is continuous in $H^1_L$.
 \item[{(iii)}] $\bar{H}$ is Gauteaux-differentiable in $H^1_L$.
 \end{enumerate}
\end{Prop}

\paragraph{\bf Proof.} The proof is similar to the proof of
Proposition~(2.3) in [13] by some obvious changes.$\square$\\

Here is the main result in this section.
\begin{Th}
 $ \bar H$ has a critical point  in $H^1_L$, that is,
there exists $0 \neq v\in H^1_L$ such that
$$\int \dDel
v.\dDel \phi dx+\int V(x)f(v)f'(v)\phi dx-\int w(x,f(v))f'(v)\phi
dx=0,$$ for every $\phi\in H^1_L$.
\end{Th}

We use the Mountain Pass Theorem (see [2], [19]) to prove Theorem
3.2.  First, let us define the Mountain Pass value,
$$C_0:=\inf_{\ga\in\Ga} \sup_{t\in[0,1]} \bar{H}(\ga(t)),$$
where
$$\Ga=\{\ga\in C([0,1], H^1_L) | \ga(0)=0, \bar {H}(\ga(1))\leq 0, \ga(1)\neq 0\}.$$
The following Lemmas are crucial for the proof of Theorem 3.2.
\begin{Lem}
 The functional $\bar{H}$ satisfies the Mountain Pass
Geometry. \end{Lem}
\paragraph{\bf Proof.}
 We need to show that there exists $0\neq v \in~H^1_L$ such that
$\bar{H}(v)\leq 0$.   Let $e\in C_{0,r}^\infty({\BR}^N)$ with
$e\not\equiv 0$ and  supp$( e)\subset\Omega$. It is easy to see
that ${H}(te)\leq 0 $ for the large values of $t.$   Consequently
  $\bar {H} (v) <0$ where $v=h(te).$ $\square$
\begin{Lem} $C_0$\; is positive.
\end{Lem}
\paragraph{\bf Proof.} Set
\begin{equation*}
S_\rho:=\{v\in H^1_L| \int |\dDel v|^2dx+\int V(x)
f(v)^2dx=\rho^2\}.
\end{equation*}

It follows from $H4$ that for a given $\epsilon >0$ there exists
$\del >0$ such that
$$G(t)\leq \f {\epsilon t^2}{2}, \quad \quad |t|\leq \del.$$

Thus
\begin{eqnarray*}
\int_{\Lambda} G(u)\,dx \leq \f {\epsilon}{2}\int_{\Lambda} u^2\,dx,
\text {  as  } \|u\|_X \leq \rho_1 \text {  with  } \rho_1 \text {
small enough.}
\end{eqnarray*}

Set $u=f(v).$ It is easy to check that $\|u\|_X\leq
\|v\|_{H^1_L({\BR}^N)}.$ Hence, it follows
\begin{eqnarray*}
\int_{\Lambda} G(f(v))\,dx \leq \f {\epsilon}{2}\int_{\Lambda}
f(v)^2\,dx, \text {  as  } \|v\|_{H^1_L({\BR}^N)} \leq \rho_1 \text
{  with  } \rho_1 \text { small enough.}
\end{eqnarray*}

Recalling that
\begin{eqnarray*}
\int_{\Lambda} f(v)^2\,dx \leq C (\int |\dDel v|^2dx+\int V(x)
f(v)^2dx),
\end{eqnarray*}
we obtain for each $v \in S_{\rho}$
\begin{eqnarray}
\int_{\Lambda} G(f(v))\,dx \leq \f {C\epsilon}{2}\rho^2,
\end{eqnarray}
for small values of $\rho.$

Also,  it follows from $(G1)$  and $(G2)$  for each $v \in S_{\rho}$
with $\rho$ small enough that
\begin{align}
\int W(x,f(v))dx &=\int_\La W(x,f(v))dx+\int_{\La^c} W(x,f(v))dx \notag\\
&\leq \int_{\Lambda} G(f(v))\,dx + \f{1}{2k}\int V(x)f(v)^2 dx \nonumber \\
& \leq \f{C\epsilon}{2}\rho^2 +\f{1}{2k}\rho^2
\end{align}
Considering  (4), (5) and the fact that $v\in S_\rho$, we obtain
\begin{align*}
\bar{H}(v)&=\f{1}{2}\int |\dDel v|^2dx+\f{1}{2} \int V(x)f(v)^2 dx-\int
W(x,f(v))dx\\
&\geq \f{1}{2}\rho^2-\f{C\epsilon}{2}\rho^2-\f{1}{2k}
\rho^2=(\f{1}{2}-\f{1}{2k})\rho^2-\f{C\epsilon}{2}\rho^2 \geq
\f{k-1}{4k}\rho^2,
\end{align*}
when $0<\rho\leq \rho_0\ll 1$ for some $\rho_0$ and $\epsilon$ small
enough. Hence, for $v\in S_\rho$ with $0<\rho\leq \rho_0$ we have
\begin{equation}
\bar{H}(v)\geq \f{k-1}{4k}\rho^2.
\end{equation}
If $\ga(1)=v$ and $\bar {H}(\ga(1))<0$ then it follows from (6) that
$$\int |\dDel v|^2dx+ \int V(x)f(v)^2dx >\rho_0^2,$$
thereby giving
$$\sup_{t\in [0,1]}\bar{H}(\ga(t)) \geq \sup_{\ga(t)\in S_{\rho_0}}
\bar{H}(\ga(t))\geq \f{k-1}{4k}\rho_0^2.$$
Therefore $C_0\geq \f{k-1}{4k}\rho_0^2>0$.$\Box$

The Mountain Pass Theorem guaranties the existence of a $(PS)_{C_0}$
sequence $\{v_n\},$ that is,  $\bar{H}(v_n)\lo C_0$ and
$\bar{H}'(v_n)\lo 0$. The following lemma states some properties of
this sequence.
\begin{Lem}Suppose $\{v_n\}$ is a $(PS)_{C_0}$ sequence. The following statements hold.
 \begin{enumerate}
 \item[{(i)}] $\{v_n\}$ is bounded in $H^1_L.$
\item[{(ii)}]  For each $\del>0$, there exists $R>4R_2$, ($R_2$ is introduced in $(A1)$ and $(A2)$) such that
$$\underset{n\rightarrow+\infty}{\mathrm{limsup}}
\int_{B_R^c} \big ( |\dDel v_n|^2+ V(x)f(v_n)^2 \big )dx<\del.$$
\item[{(iii)}]  If $v_n$ converges weakly to $v $  in $ H^1_L$, then
$$\lim_{n\rightarrow+\infty} \int w(x,f(v_n)) f(v_n)dx=\int w(x,f(v))f(v) dx.$$
\item[{(iv)}] If $v_n\geq 0$ converges  weakly to $v $  in $ H^1_L$,  then for every nonnegative test  function $\phi \in  H^1_L $ we have
$$\lim_{n\rightarrow+\infty} \braket{\bar{H}'(v_n)}{\phi}= \braket{\bar{H}'(v)}{\phi}.$$
\end{enumerate}
\end{Lem}
\paragraph{\bf Proof.} Since $\{v_n\}$ is  a $(PS)_{C_0}$ sequence, we have
\begin{eqnarray}
\bar{H}(v_n)&=&\f{1}{2} \int |\dDel v_n|^2 dx+\f{1}{2} \int V(x)f(v_n)^2
dx- \int W(x,f(v_n))dx\nonumber \\&=&C_0+o(1),
\end{eqnarray}
and
\begin{eqnarray}
\braket{\bar{H}'(v_n)}{\phi}&=&\int \dDel v_n.\dDel \phi dx+\int V(x)f(v_n)f'(v_n)\phi
dx-\int w(x, f(v_n))f'(v_n) \phi dx \nonumber \\
&=&o(\|\phi\|)
\end{eqnarray}
For part $(i),$ pick $\phi=\f{f(v_n)}{f'(v_n)}=
\sqrt{1+f(v_n)^2}f(v_n)$ as a test function. One can easily deduce
that $|\phi|_L\leq C|v_n|_L$ and
$$|\dDel \phi|=(1+\f{f(v_n)^2}{1+f(v_n)^2}) |\dDel v_n| \leq 2|\dDel v_n|,$$
which implies  $\|\phi\|\leq C\|v_n\|$. Substituting $\phi$ in (8),
gives
\begin{align}
\braket{\bar{H}'(v_n)}{
 \f{f(v_n)}{f'(v_n)}}&= \int (1+\f{f(v_n)^2}{1+f(v_n)^2}) |\dDel
v_n|^2dx+ \int V(x)f(v_n)^2dx \nonumber \\
&\quad - \int w(x,f(v_n))f(v_n)dx \nonumber \\
&=o(\|v_n\|).
\end{align}
It follows from $(G1)$ and $(G2)$ that
\begin{equation}
-\int W(x,f(v_n)) dx+\f{1}{\tht} \int w(x,f(v_n)) f(v_n) dx\geq
\f{1}{k} (\f{1}{\tht}-\f{1}{2}) \int V(x)f(v_n)^2dx
\end{equation}
 Taking into
account (7), (9) and (10), we have
\begin{align*}
C_0+o(1)+o(\|v_n\|) =&\bar{H}(v_n)-\f{1}{\tht}
\braket{\bar{H}'(v_n)}{
\f{f(v_n)}{f'(v_n)}}\\
=&\f{1}{2}\int |\dDel v_n|^2 dx+\f{1}{2} \int V(x)f(v_n)^2dx
-\int
W(x,f(v_n))dx \\
&-\f{1}{\tht} \int (1+\f{f(v_n)^2}{1+f(v_n)^2}) |\dDel v_n|^2 dx -
\f{1}{\tht}
\int V(x)f(v_n)^2dx \\
&+\f{1}{\tht} \int w(x,f(v_n))f(v_n)dx\\
=&\int(\f{1}{2}-\f{1}{\tht}(1+\f{f(v_n)^2}{1+f(v_n)^2})) |\dDel
v_n|^2dx+
(\f{1}{2}-\f{1}{\tht}) \int V(x) f(v_n)^2dx \\
&+\int (\f{1}{\tht} w(x,f(v_n))f(v_n)-W(x,f(v_n)))dx \\
\geq & (\f{1}{2}-\f{2}{\tht}) \int |\dDel v_n|^2
dx+(\f{1}{2}-\f{1}{\tht})(1-\f{1}{k}) \int V(x)f(v_n)^2dx.
\end{align*}
Since, $\f{1}{2}-\f{2}{\tht}>0$ and
$(\f{1}{2}-\f{1}{\tht})(1-\f{1}{k})>0$ it follows from the above
that   $\int |\dDel v_n|^2dx+ \int V(x) f(v_n)^2dx$ is bounded. It
proves part $(i).$

For part $(ii)$,  let $\eta_R\in C^\infty({\BR}^N,{\BR})$ be a
function satisfying $\eta_R=0$ on $B_{\f{R}{2}}$, $\eta_R=1$ on
$B_R^c$ and $|\dDel \eta_R(x)|\leq\f{C}{R}$. It follows from part
$(i)$  that $\{v_n\}$ is bounded. Hence, from (8) we have
$$ \braket{\bar{H}'(v_n)}{ \f{f(v_n)}{f'(v_n)} \eta_R}=o(1),$$
thereby giving
\begin{multline*}
\int (1+\f{f(v_n)^2}{1+f(v_n)^2}) |\dDel v_n|^2\eta_R dx+\int V(x)
f(v_n)^2\eta_R dx\\
+\int \f{f(v_n)}{f'(v_n)} \dDel v_n.\dDel \eta_R dx=\int
w(x,f(v_n))f(v_n)\eta_Rdx+o(1).
\end{multline*}
By $(G2)$, we get
$$w(x,f(v_n)) f(v_n) \leq \f{V(x)}{k} f(v_n)^2, \quad \forall x\in
B^c_{\f{R}{2}}.$$ Therefore,
\begin{align}
\int (1+\f{f(v_n)^2}{1+f(v_n)^2} )& |\dDel v_n|^2 \eta_R dx +\int
(1-\f{1}{k})V(x) f(v_n)^2\eta_R dx  \nonumber \\
&\leq \f{C}{R} \int \f{|f(v_n)|}{f'(v_n)} |\dDel v_n| dx+o(1) \nonumber  \\
&\leq \f{C}{R} \int |\dDel v_n|^2dx +\f{C}{R} \int
(|f(v_n)|^2+|f(v_n)|^4)dx+o(1).
\end{align}
Also, it follows from part $(vi)$ of Proposition 2.1 that
$\{f(v_n)\}_n$ is a bounded sequence in $L^2({\BR}^N)\cap
L^{4}({\BR}^N)$. Hence, $\int (|f(v_n)|^2+ |f(v_n)|^4) dx$ is
bounded. Therefore, it follows from (11) that
$$\underset{n\rightarrow\infty}{\mathrm{limsup}}
\int_{B_R^c} \big  (|\dDel v_n|^2dx+ V(x) f(v_n)^2 \big )dx<\del,
\quad (R>4R_2).$$ It proves part $(ii).$

 For part $(iii),$ note first that from part $(ii)$ of the present   Lemma for each $\del>0$  there exists $R>4R_2$ such that
\begin{equation}
\underset{n\rightarrow\infty}{\mathrm{limsup}}\int_{B_R^c}\big
(|\dDel v_n|^2+V(x) f(v_n)^2 \big )dx<\f{k\del}{4}.
\end{equation}
Since $B_R^c \subseteq \Lambda^c,$ it follows from $(G2)$ that
$$w(x,f(v_n)) f(v_n)\leq \f{V(x)}{k}f(v_n)^2 \quad \quad \quad \forall x \in B_R^c$$
which together with (12) imply that
\begin{equation}
\underset{n\rightarrow\infty}{\mathrm{limsup}} \int_{B_R^c}
w(x,f(v_n))f(v_n)dx\leq \f{\del}{4},
\end{equation}
 and consequently
\begin{equation*}
\int_{B_R^c} w(x,f(v))f(v)dx\leq \f{\del}{4}.
\end{equation*}
It follows from (13) and the above inequality  that
\begin{multline}
\Big |\int w(x,f(v_n))f(v_n)dx- \int w(x,f(v))f(v)dx \Big |  \\
\leq \f{\del}{2}+\Big | \int_{B_{R_1}} \big [w(x,f(v_n))f(v_n)-w(x,f(v))f(v)\big ]dx \Big |  \\
+\Big | \int_{B_R\backslash B_{R_1}}  \big
[w(x,f(v_n))f(v_n)-w(x,f(v))f(v) \big ]dx \Big |.
\end{multline}
Since $B_{R_1}\subset\La^c$, we have
$$w(x,f(v_n))f(v_n) \leq \f{V(x)}{k} f(v_n)^2, \quad \forall x\in B_{R_1}$$
Then, by the compact theorem embedding and Lebesgue Theorem, we
obtain a subsequence still denoted by $\{v_n\}$, such that
\begin{equation}
\int_{B_{R_1}} w(x,f(v_n))f(v_n)dx\lo \int_{B_{R_1}} w(x,f(v))f(v)dx.
\end{equation}

Now, we show that

\begin{equation*}
\int_{B_R\backslash \bar{B}_{R_1}} w(x,f(v_n))f(v_n)dx\lo
\int_{B_R\backslash\bar{B}_{R_1}} w(x,f(v))f(v)dx.
\end{equation*}
Since $v_n\rightharpoonup v$ weakly in $H^1_L,$ there exists a
constant $C$ such that $\|v_n\|\leq C.$ Set $u_n=f(v_n).$ An easy
computation shows that $\|u_n\|_{X}\leq \|v_n\|\leq C.$ Using
Straus's inequality (see [22]) we have

$$|u_n(x)|\leq \f{2\pi}{|x|^{\f{1}{2}}} \|u_n\|_{X}\leq \f{2\pi
C}{|x|^{\f{1}{2}}},\quad \quad \quad \quad  \forall x\neq 0,$$

from which
$$|u_n(x)|\leq \f{2\pi
C}{R_1^{\f{1}{2}}}:=\bar C,\quad \quad \quad \quad  \forall x \in
B_R\backslash\bar{B}_{R_1}.$$

From this we have
$$ |w(x,f(v_n))v_n|=|w(x,u_n)u_n|\leq \max_{x \in
B_R\backslash\bar{B}_{R_1},  t   \in [- \bar C, \bar C]} w(x,t) \bar
C:= \bar C_0 \in L^1(B_R\backslash\bar{B}_{R_1}).$$

Then, it follows from the Lebesgue dominated convergence theorem
that

\begin{equation}
\int_{B_R\backslash \bar{B}_{R_1}} w(x,f(v_n))f(v_n)dx\lo
\int_{B_R\backslash\bar{B}_{R_1}} w(x,f(v))f(v)dx.
\end{equation}

 Considering (15) and (16), it follows from (14) that
$$\underset{n\rightarrow\infty}{\mathrm{limsup}} \Big |\int w(x,f(v_n))f(v_n)dx-\int
w(x,f(v))f(v)dx \Big | \leq \f{\del}{2},$$ for every $\del>0$.
Consequently
$$\int w(x,f(v_n))f(v_n)dx\lo \int w(x,f(v))f(v)dx,$$
as $n \rightarrow \infty.$  It proves part $(iii).$

To prove  part $(iv),$ note first that $f$ is increasing and $f(0)=0$, hence $f(v_n)\geq 0$ and $ f(v)\geq0.$
For the second term on the right hand side of (8), we have
$$V(x) f(v_n) f'(v_n)\phi \leq  V(x) f(v_n)\phi,$$
and since $v_n\rightharpoonup v$ weakly in $H_1^G$, for the right
hand side of the above inequality we have
$$\lim_{n\rightarrow\infty} \int V(x) f(v_n)\phi \, dx = \int V(x)
f(v)\phi \, dx.$$  Hence by the dominated convergence theorem and
the fact that $v_n \rightarrow v$ a.e. we obtain
\begin{equation}
\lim_{n\rightarrow\infty} \int V(x) f(v_n)f'(v_n)\phi \, dx =  \int V(x) f(v)f'(v)\phi \, dx.
\end{equation}
For the third term on the right hand side of (8), we have
$$w(x,f(v_n)) f'(v_n) \phi \leq \frac{V(x)}{k}f(v_n) \phi, \quad \quad \quad \forall x \in \Lambda^c,$$
and similarly by the dominated convergence theorem, we obtain
 \begin{equation}
\lim_{n\rightarrow\infty} \int_{\Lambda^c} w(x, f(v_n))f'(v_n)\phi \, dx =  \int_{\Lambda^c} w(x, f(v))f'(v)\phi \, dx.
\end{equation}

Also, by the same argument to prove (16),  we obtain
 \begin{equation}
\lim_{n\rightarrow\infty} \int_{\Lambda} w(x, f(v_n))f'(v_n)\phi \, dx =  \int_{\Lambda} w(x, f(v))f'(v)\phi \, dx.
\end{equation}
It follows from (8) and (17)-(19) that
$$\lim_{n\rightarrow+\infty} \braket{\bar{H}'(v_n)}{\phi}=\braket{\bar{H}'(v)}{\phi}.$$
It proves part $(iv)$. $\square$
\begin{Lem} If $\{v_n\}$ is a $(PS)_{C_0}$ sequence, then $v_n$
converges to $v\in H^1_L$. Consequently
$\bar{H}(v)=\lim_{n\rightarrow+\infty} \bar{H}(v_n)$ and
$\bar{H}'(v)=0$. \end{Lem}
\paragraph{\bf Proof.}  It follows from part $(i)$ of Lemma 3.5 that $v_n$ is a bounded sequence in  $H^1_L.$
Hence, there exists $v \in  H^1_L $ such that, up to a subsequence,
$v_n\rightharpoonup v$ weakly in $H^1_L$ and  $v_n\rightarrow v$
a.e. in ${\BR}^N.$ Since we may replace $v_n$ by $|v_n|,$ we assume
$v_n\geq  0$ and $v\geq 0.$  Since, $\{v_n\}$ is a $(PS)_{C_0}$
sequence we have
\begin{align}
o(\|v_n\|)&= \braket{\bar{H}'(v_n)}{
 \f{f(v_n)}{f'(v_n)}}\\&= \int (1+\f{f(v_n)^2}{1+f(v_n)^2}) |\dDel
v_n|^2dx+ \int V(x)f(v_n)^2dx \quad - \int w(x,f(v_n))f(v_n)dx
\nonumber
\end{align}
and
 \begin{align}
o(\|v\|)=\braket{\bar{H}'(v_n)}{
 \f{f(v)}{f'(v)}}.
\end{align}
It follows from  part $(iv)$ of Lemma 3.5  and (21) that
 \begin{align}
\braket{\bar{H}'(v_n)}{
 \f{f(v)}{f'(v)}} =& \braket{\bar{H}'(v)}{
 \f{f(v)}{f'(v)}}+o(\|v\|) \nonumber \\=&\int (1+\f{f(v)^2}{1+f(v)^2}) |\dDel
v|^2dx+ \int V(x)f(v)^2dx \nonumber \\
& - \int w(x,f(v))f(v)dx +o(\|v\|)
\end{align}
In this step, we show that
\begin{align*}
 \int \frac{f(v)^2 |\dDel v|^2}{1+f(v)^2} \, dx \leq \liminf_{n\rightarrow \infty } \int \frac{f(v_n)^2 |\dDel v_n|^2}{1+f(v_n)^2} \, dx.
\end{align*}
 Set $u_n=f(v_n)$ and $u=f(v).$ A direct
computation shows that
$$\int |\dDel u_n^2|^2 \, dx=4 \int \frac{f(v_n)^2 |\dDel v_n|^2}{1+f(v_n)^2} \, dx\leq 4\|v_n\|^2.$$
Also, from part (vi) of Proposition 2.1 we have

$$\int u_n^4 \, dx = \int f(v_n)^4 \, dx \leq C \|v_n\|^4.$$
 Set $w_n=u_n^2.$ It follows from the above that $\{w_n\}_n$ is a
bounded sequence in $H^1({\BR}^N).$ Hence, up to a subsequence $w_n
\rightharpoonup w$  weakly in $H^1({\BR}^N)$ and $w_n \rightarrow w$
a.e. in ${\BR}^N.$ It follows $w=u^2.$ Also, by the lower semi
continuity of the  norm in $H^1({\BR}^N),$ we obtain
$$\int |\dDel w|^2 \, dx \leq \liminf_{n\rightarrow \infty } \int |\dDel w_n|^2 \, dx.$$
Plug  $w_n=u_n^2$ and $w=u^2$ in this inequality to get
$$\int |\dDel u^2 |^2 \, dx \leq \liminf_{n\rightarrow \infty } \int |\dDel u_n^2 |^2 \, dx.$$
Substituting $u_n=f(v_n)$ and $u=f(v)$ in the above  inequality
gives
\begin{align}
 \int \frac{f(v)^2 |\dDel v|^2}{1+f(v)^2} \, dx \leq \liminf_{n\rightarrow \infty } \int \frac{f(v_n)^2 |\dDel v_n|^2}{1+f(v_n)^2} \, dx.
\end{align}
Also, lower  semi continuity and Fatou's Lemma imply
\begin{align}
\int |\dDel v|^2dx &\leq \underset{n\rightarrow\infty}{\mathrm{liminf}} \int
|\dDel v_n|^2dx, \\
\int V(x)L(v)dx &\leq
\underset{n\rightarrow\infty}{\mathrm{liminf}}\int V(x)L(v_n)dx.
\end{align}
Up to a subsequence  one can assume
\begin{align}
\liminf_{n\rightarrow \infty }\int
|\dDel v_n|^2dx&= \lim_{n\rightarrow \infty } \int
|\dDel v_n|^2dx \\
\liminf_{n\rightarrow \infty }\int V(x)L(v_n)dx&= \lim_{n\rightarrow
\infty }\int
V(x)L(v_n)dx.\\
\liminf_{n\rightarrow \infty } \int \frac{f(v_n)^2 |\dDel v_n|^2}{1+f(v_n)^2} \, dx&=\lim_{n\rightarrow \infty } \int \frac{f(v_n)^2 |\dDel v_n|^2}{1+f(v_n)^2} \, dx.
\end{align}
It follows from (23)-(28) that there exist nonnegative numbers
$\delta_1, \delta_2$ and $\delta_3$ such that
\begin{align}
\lim_{n\rightarrow \infty }\int
|\dDel v_n|^2dx&= \int
|\dDel v|^2dx+\delta_1 \\
\lim_{n\rightarrow \infty }\int V(x)L(v_n)dx&= \int
V(x)L(v)dx+\delta_2.\\
\lim_{n\rightarrow \infty } \int \frac{f(v_n)^2 |\dDel v_n|^2}{1+f(v_n)^2} \, dx&= \int \frac{f(v)^2 |\dDel v|^2}{1+f(v)^2} \, dx+\delta_3.
\end{align}
Now, we show that $\delta_1=\delta_2=\delta_3=0.$  It follows from
part $(iii)$ of Lemma 3.5 that
$$\int w(x,f(v_n)) f(v_n)dx\lo \int w(x,f(v))f(v)dx.$$
which together with (20) and (22) imply
 \begin{align*}
 \lim_{n\rightarrow \infty } \Big\{\int (1+\f{f(v_n)^2}{1+f(v_n)^2}) |\dDel
v_n|^2dx\\+ \int V(x)f(v_n)^2dx \Big \}&=
\quad   \lim_{n\rightarrow \infty } \int w(x,f(v_n))f(v_n)dx \\
&= \int w(x,f(v))f(v)dx\\
&=\int (1+\f{f(v)^2}{1+f(v)^2}) |\dDel
v|^2dx+ \int V(x)f(v)^2dx
\end{align*}
Taking into account (29), (30) and (31) the above limit implies
$\delta_1=\delta_2=\delta_3=0.$ Therefore, it follows from (29) and
(30) that
\begin{align*}
\int |\dDel v|^2dx&=\underset{n\rightarrow\infty}{\mathrm{lim}} \int
|\dDel
v_n|^2dx \\
\int V(x)L(v)dx&=\underset{n\rightarrow\infty}{\mathrm{lim}}\int
V(x)L(v_n)dx.
\end{align*}
By Proposition 2.1, $v_n\lo v$ in $E_L$ and we have $\dDel v_n\lo
\dDel v$ in $L^2$. Hence $v_n\lo v$ in $H^1_L$. $\square$
\paragraph{\bf Proof of Theorem 3.2.} The proof is a direct
consequence of Lemmas 3.3, 3.4 and 3.5. $\square$
\section{Proof of Theorem 1.1}
To prove Theorem 1.1, note first that every critical point of the
functional  $\bar{J}_\eps$ corresponds to a weak solution of problem
(2). Thus, we need to find a critical point for the functional
$\bar{J}_\eps.$ To do this, we shall show that the functionals
$\bar{J}_\eps$ and $\bar{H}_\eps$ will coincide for the small values
of $\eps$. Hence, every critical point of $\bar{H}_\eps$ will be a
critical point of $\bar{J}_\eps.$ Also, it follows from Theorem 3.2
that $\bar{H}_\eps$ has a
nontrivial critical point for every $\eps>0$. \\

Without loss of generality,  we may assume $\eps^2$ instead of
$\eps$ in the functionals $\bar{H}_\eps$ and $\bar{J}_\eps$, i.e.

$$\bar{H}_\eps (v)=\f{\eps^2}{2} \int |\dDel v|^2+\f {1}{2}\int
V(x)f(v)^2dx-\int W(x,f(v))dx,$$
and

$$\bar {J}_\eps(v)=\f{\eps^2}{2} \int_{{\BR}^N} |\dDel v|^2dx+ \f{1}{2} \int_{{\BR}^N}V(x)
f(v)^2 dx- \int_{{\BR}^N} G(f(v))dx.$$

It follows from Theorem 3.2 that    there exists a critical point
$v_\eps\in H^1_L$ of $\bar{H}_\eps(v)$ for each $\eps>0$. Set
$u_\eps=f(v_\eps)$.

The following Lemmas are crucial for the proof of Theorem 1.1.

\begin{Lem}
 The sequence $\{u_\eps\}_{\eps>0}$ is strongly convergent
to $0$ when $\eps\lo 0$, in $H^1({\BR}^N)$, i.e.
$$\|u_\eps\|_{H^1}\lo 0 \quad \text{as}\quad \eps\lo 0.$$
\end{Lem}
\paragraph{\bf Proof.} Let $0\not\equiv \phi\in C_{0,r}^\infty({\BR}^N)$ be a
non-negative function with supp$(\phi)\subset\Om$ and $H_1(\phi)\leq
0.$  Set $\ga_1 (t):= h(t \phi).$ Hence, we have $$ \bar { H}_\eps
(\gamma_1 (1))= \bar {H}_\eps (h (\phi))= H_\eps(\phi)\leq
H_1(\phi)\leq0.$$

It follows from the definition of the Mountain Pass value that
$$\bar{H}_\eps(v_\eps)=\inf_{\ga\in\Ga} \sup_{t\in[0,1]}
\bar{H_\eps}(\ga(t))\leq \sup_{t\in[0,1]}
\bar{H_\eps}(\ga_1(t))=\sup_{t\in[0,1]} \bar{H_\eps}(h(t
\phi))=\sup_{t\in[0,1]} {H_\eps}(t\phi).$$

Therefore, we obtain
\begin{align}
\bar{H}_\eps(v_\eps)  &\leq
\sup_{t\in[0,1]}H_\eps(t\phi) \nonumber \\
&=\sup_{t\in [0,1]} \f{\eps^2 t^2}{2} \int |\dDel \phi|^2+\f{\eps^2
t^4}{2} \int
|\phi|^2 |\dDel\phi|^2- \int G(t\phi)dx  \nonumber\\
&\leq \sup_{t\in[0,1]}\f{\eps^2 t^2}{2} \int (1+|\phi|^2)
|\dDel\phi|^2dx -
\int G(t\phi)dx \nonumber \\
&=\f{\eps^2 t_{\eps}^2}{2} \int (1+|\phi|^2) |\dDel\phi|^2dx - \int
G(t_{\eps}\phi)dx
\end{align}
for some $0<t_{\eps}<1.$  On the other hand we have
\begin{eqnarray*}
\eps^2 \int (1+|\phi|^2) |\dDel\phi|^2dx = \int \f
{g(t_{\eps}\phi)\phi}{t_{\eps}} dx.
\end{eqnarray*}
Choosing $\Omega_1 \subseteq \Omega$ such that $\phi(x)\geq
\phi_0>0,  \forall x \in \Omega_0,$ it follows

\begin{eqnarray*}
\eps^2 \int (1+|\phi|^2) |\dDel\phi|^2dx \geq \int_{\Omega_0} \f
{g(t_{\eps}\phi)\phi}{t_{\eps}} dx \geq \phi_0^2 \int_{\Omega_0}
\f{g(t_{\eps}\phi)}{t_{\eps}\phi} dx.
\end{eqnarray*}

Thus,  from the above inequalities and Conditions $H1-H3$ we obtain
$t_{\eps}\rightarrow 0$ as $\eps \rightarrow 0.$
 Now, as in the
proof of part $(i)$ of Lemma 3.5  we obtain
\begin{align}
\bar{H}_\eps(v_\eps) & = \bar{H}_\eps(v_\eps)-\f{1}{\tht} \braket {
\bar{H}'(v_\eps)}{v_\eps} \nonumber \\
&\geq \eps^2(\f{1}{2}-\f{2}{\tht}) \int |\dDel v_n|^2dx+
(\f{1}{2}-\f{1}{\tht})(1-\f{1}{k}) \int V(x)f(v_n)^2dx.
\end{align}
Combining (32) and (33), we get
\begin{eqnarray*}
\eps^2(\f{1}{2}-\f{2}{\tht}) \int |\dDel v_n|^2 dx+
(\f{1}{2}-\f{1}{\tht}) (1-\f{1}{k}) \int V(x) |f(v_n)|^2 dx &\leq &
\f{\eps^2 t_{\eps}^2}{2} \int (1+|\phi|^2) |\dDel\phi|^2dx \\
& &-\int
G(t_{\eps}\phi)dx\\
&\leq &  \f{\eps^2 t_{\eps}^2}{2} \int (1+|\phi|^2) |\dDel\phi|^2dx.
\end{eqnarray*}
Therefore
\begin{equation}
(\f{1}{2}-\f{2}{\tht}) \int |\dDel v_n|^2 dx+ (\f{1}{2}-\f{1}{\tht})
(1-\f{1}{k}) \int V(x) f(v_n)^2 dx  \leq \f{t_{\eps}^2}{2} \int
(1+|\phi|^2) |\dDel\phi|^2dx.
\end{equation}
Hence, substituting $u_\eps=f(v_\eps)$ in (34) implies
$$\int (1+|u_\eps|^2) |\dDel u_\eps|^2dx+\int V(x) |u_\eps|^2dx\leq
\f{t_{\eps}^2}{2} \int (1+|\phi|^2) |\dDel\phi|^2dx.$$ Therefore
$$\|u_\eps\|_{H^1}\lo 0 \quad \text{as}\quad \eps\lo 0.$$
$\square$\\

\begin{Lem}
 For every compact set $Q\subset{\BR}^N$ such that
$0\not\in Q$, $\|u_\eps\|_{L^\infty(Q)}\lo 0$ as $\eps\lo 0$.
\end{Lem}
\paragraph{\bf Proof.} For each $\eps>0$, it follows from Straus's inequality  that
$$0\leq u_\eps(x)\leq \f{2 \pi}{|x|^{\f{1}{2}}} \|u_\eps\|_{H^1({\BR}^N)}\quad
\forall x\neq 0,$$ which together with the result of Lemma 4.1
obviously means
$$\|u_\eps\|_{L^\infty(Q)}\lo 0 \quad \text{as}\quad \eps\lo 0.$$ $\square$
\paragraph{\bf Proof of Theorem 1.1.} By Lemma 4.2 we have
\begin{equation}
M_\eps:=\max_{x\in \bar \La} f(v_\eps)\lo 0 \quad \text{as}\quad
\eps\lo 0.
\end{equation}

From (35)  there exists $\eps_0>0$ such that $\max_{x\in \bar
\La}f(v_\eps)<a$
 for every $0<\eps<\eps_0$. Using the
test function $\phi=\f{(f(v_\eps)-a)_+}{f'(v_\eps)}$, we get
\begin{eqnarray*}
0=\braket{\bar{H}'_\eps(v_\eps)}{\phi} = \int_{\textit{F}}
\eps^2(1+\f{f(v_\eps)^2}{1+f(v_\eps)^2}) |\dDel v_\eps|^2 &+&
\int_{{\BR}^N \backslash \bar \La} V(x)f(v_\eps)(f(v_\eps)-a)_+dx\\
&-&\int_{{\BR}^N \backslash \bar \La} w(x,f(v_\eps))
(f(v_\eps)-a)_+dx
\end{eqnarray*}
where $\textit{F}=({{\BR}}^N \backslash \bar \La)\cap \{x|
f(v_\eps)\geq a\}$. From $(G2)$, we have
$$V(x)f(v_\eps)(f(v_\eps)-a)_+- w(x,f(v_\eps))(f(v_\eps)-a)_+\geq 0,
\quad \forall x\in\La^c.$$
Thus,
$$\eps^2 \int_{\textit{F}} (1+\f{f(v_\eps)^2}{1+f(v_\eps)^2}) |\dDel v_\eps|^2dx=0,$$
from which we obtain
$$f(v_\eps)\leq a, \quad \forall x\in {\BR}^N \backslash \bar \La.$$
Therefore
$$w(x,f(v_\eps))= g(f(v_\eps)) , \quad \forall x\in {\BR}^N \backslash \bar \La,$$
and we conclude that
$$\eps^2 \int \dDel v_\eps. \dDel \xi dx+ \int V(x) f(v_\eps) f'(v_\eps) \xi dx =
\int g(f(v_\eps))  f'(v_\eps) \xi dx$$ for every $ \xi\in H^1_L$ and
$\eps\in (0,{\eps}_0)$. Therefore, $\bar {J}_\eps(v)$ has a critical
point $v_\eps$ in $H^1_L$  for every $\eps\in (0,{\eps}_0)$.
$\square$%

\end{document}